\documentclass[12pt,a4paper]{article}
\usepackage{amsmath}
\usepackage{amssymb}
\usepackage{t1enc}
\usepackage[latin1]{inputenc}
\usepackage[english]{babel}
\usepackage{amsfonts}
\usepackage{latexsym}
\usepackage{bm}
\newtheorem{theorem}{Theorem}[section]

\newtheorem{lemma}[theorem]{Lemma}
\newtheorem{cor}[theorem]{Corollary}

\begin{document}
\title{A cross-intersection theorem for subsets of a set}

\author{Peter Borg\\[5mm]
Department of Mathematics, University of Malta, Malta\\
\texttt{peter.borg@um.edu.mt}}
\date{}
\maketitle

\begin{abstract}
Two families $\mathcal{A}$ and $\mathcal{B}$ of sets are said to be \emph{cross-intersecting} if each member of $\mathcal{A}$ intersects each member of $\mathcal{B}$. For any two integers $n$ and $k$ with $0 \leq k \leq n$, let ${[n] \choose \leq k}$ denote the family of all subsets of $\{1, \dots, n\}$ of size at most $k$. 
We show that if $\mathcal{A} \subseteq {[m] \choose \leq r}$, $\mathcal{B} \subseteq {[n] \choose \leq s}$, and $\mathcal{A}$ and $\mathcal{B}$ are cross-intersecting, then 
%
\[|\mathcal{A}||\mathcal{B}| \leq \sum_{i=0}^r {m-1 \choose i-1} \sum_{j=0}^s {n-1 \choose j-1},\] 
and equality holds if $\mathcal{A} = \{A \in {[m] \choose \leq r} \colon 1 \in A\}$ and $\mathcal{B} = \{B \in {[n] \choose \leq s} \colon 1 \in B\}$. Also, we generalise this to any number of such cross-intersecting families.
\end{abstract}

\section{Basic definitions and notation} \label{Intro}
Unless otherwise stated, we shall use small letters such as $x$ to
denote elements of a set or non-negative integers or functions,
capital letters such as $X$ to denote sets, and calligraphic
letters such as $\mathcal{F}$ to denote \emph{families}
(i.e.~sets whose members are sets themselves). It is to be
assumed that arbitrary sets and families are \emph{finite}. We
call a set $A$ an \emph{$r$-element set} 
if its size $|A|$ is $r$, i.e.~if it contains exactly $r$ elements (also called members).

$\mathbb{N}$ denotes the set $\{1, 2, \dots\}$ of all positive
integers. For any integer $n \geq 0$, the set $\{i \in \mathbb{N}
\colon i \leq n\}$ is denoted by $[n]$; note that $[0]$ is the empty set $\emptyset$. For a set $X$, the \emph{power set of $X$} (i.e.~$\{A \colon A \subseteq X\}$) is denoted by $2^X$. For any integer $r \geq 0$, the family of all $r$-element subsets of $X$ is denoted by $X \choose r$, and the family of all subsets of $X$ of size at most $r$ is denoted by $X \choose \leq r$. So ${X \choose \leq r} = \{A \subseteq X \colon |A| \leq r\} = \bigcup_{i = 0}^r {X \choose i}$.

We will denote the union of all sets in a family $\mathcal{F}$ (i.e.~$\bigcup_{F \in \mathcal{F}} F$) by $U(\mathcal{F})$. If $x$ is an element of a set $X$, then we
denote the family of those sets in $\mathcal{F}$ which contain
$x$ (i.e.~$\{F \in \mathcal{F} \colon x \in F\}$) by $\mathcal{F}(x)$, and we call $\mathcal{F}(x)$ a \emph{star of $\mathcal{F}$}. So $\mathcal{F}(x)$ is the empty set $\emptyset$ if and only if $x$ is not in $U(\mathcal{F})$.

We say that a set $A$ \emph{intersects} a set $B$ if $A$ and $B$ contain at least one common element (i.e.~$A \cap B \neq \emptyset$). 

A family $\mathcal{A}$ is said to be \emph{intersecting} if every two members of $\mathcal{A}$ intersect. Note that the stars of a family $\mathcal{F}$ are the simplest intersecting subfamilies of $\mathcal{F}$.

If $\mathcal{A}$ and $\mathcal{B}$ are families such that each set in $\mathcal{A}$ intersects each set in $\mathcal{B}$, then $\mathcal{A}$ and $\mathcal{B}$ are said to be \emph{cross-intersecting}. In general, families $\mathcal{A}_1, \dots, \mathcal{A}_k$ are said to be \emph{cross-intersecting} if for every $i, j \in [k]$ with $i \neq j$, each set in $\mathcal{A}_i$ intersects each set in $\mathcal{A}_j$. Note that if $\mathcal{A}_1 = \dots = \mathcal{A}_k = \mathcal{S}$ for some star $\mathcal{S}$ of a family $\mathcal{F}$, then $\mathcal{A}_1, \dots, \mathcal{A}_k$ are cross-intersecting.

One of the most popular endeavours in extremal set theory is that of determining the size of a largest intersecting subfamily of a given family $\mathcal{F}$. This started in \cite{EKR}, which features the classical result, known as the Erd\H os-Ko-Rado (EKR) Theorem, that says that if $r \leq n/2$, then the size of a largest intersecting subfamily of ${[n] \choose r}$ is the size ${n-1 \choose r-1}$ of any star of ${[n] \choose r}$. 
There are various proofs of the EKR Theorem, two of which are particularly short and beautiful: Katona's \cite{K}, introducing the elegant cycle method, and Daykin's \cite{D}, using the fundamental Kruskal-Katona Theorem \cite{Ka,Kr}. The EKR Theorem gave rise to some of the highlights in extremal set theory \cite{Kat,HM,F_t1,W,AK1} and inspired many results, including generalisations \cite{T,Borg1}, that establish how large a system of sets can be under certain intersection conditions; see \cite{DF,F,F2,HT,Borg7}.

For intersecting subfamilies of a given family $\mathcal{F}$,
the natural question to ask is how large they can be. For
cross-intersecting families, two natural parameters arise: the
sum and the product of sizes of the families (note that the product of sizes of $k$ families $\mathcal{A}_1, \dots, \mathcal{A}_k$ is the number of $k$-tuples $(A_1, \dots, A_k)$ such that $A_i \in \mathcal{A}_i$ for each $i \in [k]$). It is therefore natural to consider the problem of maximising the sum or the product of sizes of $k$ cross-intersecting subfamilies (not necessarily distinct or non-empty) of a given family $\mathcal{F}$. In \cite{Borg8} this problem is analysed in a general way, and it is shown that for $k$ sufficiently large it reduces to the intersection problem (i.e.~the problem of maximising the size of an intersecting subfamily of $\mathcal{F}$). 
This problem has recently attracted much attention. Solutions have been obtained for various families (see \cite{Borg8}), including ${[n] \choose r}$ \cite{H,Pyber,MT,Bey3,Borg4,WZ}, $2^{[n]}$ \cite{MT2,Borg8},  
and families of vector spaces \cite{ST}. 

In this paper we address the maximum product problem for ${[n] \choose \leq r}$. We will actually solve the more general problem where the cross-intersecting families do not necessarily come from the same family of this kind. The following is our main result.


\begin{theorem}\label{mainthm} If $m, n \in \mathbb{N}$, $r \in [m]$, $s \in [n]$, $\mathcal{A} \subseteq {[m] \choose \leq r}$, $\mathcal{B} \subseteq {[n] \choose \leq s}$, and $\mathcal{A}$ and $\mathcal{B}$ are cross-intersecting, then 
%
\[|\mathcal{A}||\mathcal{B}| \leq \sum_{i=0}^r {m-1 \choose i-1} \sum_{j=0}^s {n-1 \choose j-1},\] 
and equality holds if $\mathcal{A} = \{A \in {[m] \choose \leq r} \colon 1 \in A\}$ and $\mathcal{B} = \{B \in {[n] \choose \leq s} \colon 1 \in B\}$.
\end{theorem}
The proof is given in Section~\ref{Proof1}.  Theorem~\ref{mainthm} yields the following generalisation for $k \geq 2$ cross-intersecting families.

\begin{theorem} \label{maincor2} Let $n_1, \dots, n_k \in \mathbb{N}$, $k \geq 2$. For each $i \in [k]$, let $r_i \in [n_i]$ and $\mathcal{A}_i \subseteq {[n_i] \choose \leq r_i}$.  If $\mathcal{A}_1, \dots, \mathcal{A}_k$ are cross-intersecting, then
\[\prod_{i=1}^k |\mathcal{A}_i| \leq \prod_{i=1}^k \left( \sum_{j = 1}^{r_i} {n_i-1 \choose j-1} \right),\]
and equality holds if $\mathcal{A}_i = \{A \in {[n_i] \choose \leq r_i} \colon 1 \in A\}$ for each $i \in [k]$.
\end{theorem}
\textbf{Proof.} For each $i \in [k]$, let $a_i =
|\mathcal{A}_i|$ and let $s_i = |\{A \in {[n_i] \choose \leq r_i} \colon 1 \in A\}|$. By Theorem~\ref{mainthm}, $a_ia_j \leq s_is_j$ for any $i, j \in [k]$ with $i \neq j$. Let mod$^*$ represent the usual \emph{modulo operation} with the exception that for every two integers $x$ and $y > 0$, $(xy) \, {\rm mod}^* \, y$ is $y$ rather than $0$. We have
\begin{align} \left( \prod_{i=1}^k a_i \right)^2 &=
(a_1a_2)(a_{3 \, {\rm mod}^* \, k}a_{4 \, {\rm mod}^* \,
k})\cdots(a_{(2k-1) \, {\rm mod}^* \, k} a_{(2k) \, {\rm mod}^*
\, k}) \nonumber \\
&\leq (s_1s_2)(s_{3 \, {\rm mod}^* \, k}s_{4 \, {\rm mod}^* \,
k})\cdots(s_{(2k-1) \, {\rm mod}^* \, k} s_{(2k) \, {\rm mod}^*
\, k}) = \left( \prod_{i=1}^k s_i \right)^2 \nonumber
\end{align}
So $\prod_{i=1}^k a_i \leq \prod_{i=1}^k s_i$. Hence the result.~\hfill{$\Box$} \\


Since ${[n] \choose \leq n} = 2^{[n]}$, we have the following.

\begin{cor} \label{maincor3} If $n_1, \dots, n_k \in \mathbb{N}$, $k \geq 2$, $\mathcal{A}_i \subseteq 2^{[n_i]}$ for each $i \in [k]$, and $\mathcal{A}_1, \dots, \mathcal{A}_k$ are cross-intersecting, then
\[\prod_{i=1}^k |\mathcal{A}_i| \leq \prod_{i=1}^k 2^{n_i-1} = 2^{\left(\sum_{i=1}^k n_i \right) - k},\]
and equality holds if $\mathcal{A}_i = \left\{A \in 2^{[n_i]} \colon 1 \in A\right\}$ for each $i \in [k]$.
\end{cor}
%

We now start working towards the proof of Theorem~\ref{mainthm}. Our approach is based on the idea of generalising the setting enough for induction to work, and we will use the \emph{compression} technique (see Section~\ref{Compsection}) together with a new alteration method. Basically, the approach is as follows. We use induction on $m + n$. The challenging part is the case $m = n$. The first problem that arises is that we can have a set $A \in \mathcal{A}$ and a set $B \in \mathcal{B}$ that intersect only in $n$; in this case, we cannot simply remove $n$ and apply the induction hypothesis. Thus we consider two alterations: removing $A$ from $\mathcal{A}$ and adding $B \backslash \{n\}$ to $\mathcal{B}$, and removing $B$ from $\mathcal{B}$ and adding $A \backslash \{n\}$ to $\mathcal{A}$. This yields two new pairs of cross-intersecting families. The second problem is that the product of the sizes of a new pair obtained in this way may become smaller. However, upon subdividing the problem appropriately and applying, for each resulting case, the right inequalities that arise from the alterations, we manage to overcome this difficulty. 

The next section is dedicated to some basic results, used in the proof of Theorem~\ref{mainthm}, about the compression operation.

\section{The compression operation}
\label{Compsection}
For any $i, j \in [n]$, let $\delta_{i,j}
\colon 2^{[n]} \rightarrow 2^{[n]}$ be defined by
\[ \delta_{i,j}(A) = \left\{ \begin{array}{ll}
(A \backslash \{j\}) \cup \{i\} & \mbox{if $j \in A$ and $i \notin
A$};\\
A & \mbox{otherwise,}
\end{array} \right. \]
and let $\Delta_{i,j} \colon 2^{2^{[n]}} \rightarrow 2^{2^{[n]}}$ be the compression operation defined by
\[\Delta_{i,j}(\mathcal{A}) = \{\delta_{i,j}(A) \colon A \in
\mathcal{A}, \delta_{i,j}(A) \notin \mathcal{A}\} \cup \{A \in
\mathcal{A} \colon \delta_{i,j}(A) \in \mathcal{A}\}.\]

This operation was introduced in the original proof \cite{EKR} of the EKR Theorem and is a very useful tool in extremal set theory. A survey on the properties and uses of compression (also called \emph{shifting}) operations is given in \cite{F}; \cite{HST} is also recommended. 

Note that $|\Delta_{i,j}(\mathcal{A})| = |\mathcal{A}|$. We will need the following basic result, which we prove for completeness.

\begin{lemma}\label{compcross} Let $\mathcal{A}$ and $\mathcal{B}$ be cross-intersecting subfamilies of $2^{[n]}$, and let $i, j \in [n]$. Then $\Delta_{i,j}(\mathcal{A})$ and $\Delta_{i,j}(\mathcal{B})$ are cross-intersecting subfamilies of $2^{[n]}$.
\end{lemma}
\textbf{Proof.} 
Suppose $A \in \Delta_{i,j}(\mathcal{A})$ and $B \in \Delta_{i,j}(\mathcal{B})$. If $A \in \mathcal{A}$ and $B \in \mathcal{B}$, then $A \cap B \neq \emptyset$ since $\mathcal{A}$ and $\mathcal{B}$ are cross-intersecting. Suppose $A \notin \mathcal{A}$ or $B \notin \mathcal{B}$; we may assume that $A \notin \mathcal{A}$. Then $A = \delta_{i,j}(A') \neq A'$ for some $A' \in \mathcal{A}$.  So $i \notin A'$, $j \in A'$, $i \in A$ and $j \notin A$. Suppose $A \cap B = \emptyset$. So $i \notin B$ and hence $B \in \mathcal{B}$. So $B \in \mathcal{B} \cap \Delta_{i,j}(\mathcal{B})$ and hence $B, \delta_{i,j}(B) \in \mathcal{B}$. So $A' \cap B \neq \emptyset$ and $A' \cap \delta_{i,j}(B) \neq \emptyset$. From $A \cap B = \emptyset$ and $A' \cap B \neq \emptyset$ we get $A' \cap B = \{j\}$, but this yields the contradiction that $A' \cap \delta_{i,j}(B) = \emptyset$.~\hfill{$\Box$}\\


If $i < j$, then we call $\Delta_{i,j}$ a \emph{left-compression}. We say that a family $\mathcal{F} \subseteq 2^{[n]}$ is \emph{compressed} if $\Delta_{i,j}(\mathcal{F}) = \mathcal{F}$ for every $i,j \in [n]$ with $i < j$ (i.e.~if $\mathcal{F}$ is invariant under left-compressions). Therefore, a family $\mathcal{F} \subseteq 2^{[n]}$ is compressed if and only if $(F \backslash \{j\}) \cup \{i\} \in \mathcal{F}$ for every $i, j \in [n]$ and every $F \in \mathcal{F}$ such that $i < j \in F$ and $i \notin F$. The families $2^{[n]}$, ${[n] \choose r}$ and ${[n] \choose \leq r}$ are important examples of compressed families. 

Suppose that a subfamily $\mathcal{A}$ of $2^{[n]}$ is not compressed. Then $\mathcal{A}$ can be transformed to a compressed family through left-compressions as follows. Since $\mathcal{A}$ is not compressed, we can find a left-compression that changes $\mathcal{A}$, and we apply it to $\mathcal{A}$ to obtain a new subfamily of $2^{[n]}$. We keep on repeating this (always applying a left-compression to the last family obtained) until we obtain a subfamily of $2^{[n]}$ that is invariant under every left-compression (such a point is indeed reached, because if $\Delta_{i,j}(\mathcal{F}) \neq \mathcal{F} \subseteq 2^{[n]}$ and $i < j$, then $0 < \sum_{G \in \Delta_{i,j}(\mathcal{F})} \sum_{b \in G} b < \sum_{F \in \mathcal{F}} \sum_{a \in F} a$).

Now consider $\mathcal{A}, \mathcal{B} \subseteq 2^{[n]}$ such that $\mathcal{A}$ and $\mathcal{B}$ are cross-intersecting. Then, by Lemma~\ref{compcross}, we can obtain $\mathcal{A}^*, \mathcal{B}^* \subseteq 2^{[n]}$ such that $\mathcal{A}^*$ and $\mathcal{B}^*$ are compressed and cross-intersecting, $|\mathcal{A}^*| = |\mathcal{A}|$ and $|\mathcal{B}^*| = |\mathcal{B}|$. Indeed, similarly to the above procedure, if we can find a left-compression that changes at least one of $\mathcal{A}$ and $\mathcal{B}$, then we apply it to both $\mathcal{A}$ and $\mathcal{B}$, and we keep on repeating this (always performing this on the last two families obtained) until we obtain $\mathcal{A}^*, \mathcal{B}^* \subseteq 2^{[n]}$ such that both $\mathcal{A}^*$ and $\mathcal{B}^*$ are invariant under every left-compression.

\section{Proof of Theorem~\ref{main}} \label{Proof1}

We now prove Theorem~\ref{main}. We actually prove a more general result, Theorem~\ref{main} below, for which we need the following definitions and notation.

A family $\mathcal{H}$ is said to be a \emph{hereditary family}
(also called an \emph{ideal}, a \emph{downset}, and an \emph{abstract simplicial complex}) if for each $H \in \mathcal{H}$, all the subsets of $H$ are members of $\mathcal{H}$. 
Clearly, a family is hereditary if and only if it is a union of power sets. For a family $\mathcal{F}$, a \emph{base of $\mathcal{F}$} is a member of $\mathcal{F}$ that is not a subset of another member of $\mathcal{F}$. So a hereditary family is the union of power sets of its bases. 

Recall the notation $\mathcal{F}(x)$, introduced in Section~\ref{Intro}. We will prove the following result.


\begin{theorem}\label{main} If $m, n \in \mathbb{N}$, $\mathcal{A} \subseteq \mathcal{G} \subseteq 2^{[m]}$, $\mathcal{B} \subseteq \mathcal{H} \subseteq 2^{[n]}$, $\mathcal{G}$ and $\mathcal{H}$ are hereditary and compressed, and $\mathcal{A}$ and $\mathcal{B}$ are cross-intersecting,
then
\[|\mathcal{A}||\mathcal{B}| \leq \left| \mathcal{G}(1) \right| \left| \mathcal{H}(1) \right|,\]
and equality holds if $\mathcal{A} = \mathcal{G}(1)$ and $\mathcal{B} = \mathcal{H}(1)$.
\end{theorem}
This generalises Theorem~\ref{mainthm} because ${[m] \choose \leq r}$ is a compressed hereditary subfamily of $2^{[m]}$ and ${[n] \choose \leq s}$ is a compressed hereditary subfamily of $2^{[n]}$.


%
%

We start by establishing the following property of hereditary families, which will have a key role.

\begin{lemma}\label{her-lemma} If $X$ and $Y$ are bases of a hereditary family $\mathcal{H}$, $x \in X$ and $x \notin Y$, then
\[|\mathcal{H}(x)| < \frac{1}{2}|\mathcal{H}|.\]
\end{lemma}
\textbf{Proof.} Let $\mathcal{I} = \mathcal{H} \backslash
\mathcal{H}(x)$. So $\mathcal{I} = \{H \in \mathcal{H} \colon x
\notin H\}$. Since $\mathcal{H}$ is hereditary, $A \backslash
\{x\} \in \mathcal{I}$ for any $A \in \mathcal{H}(x)$. Thus, we
can define a function $f \colon \mathcal{H}(x) \rightarrow
\mathcal{I}$ by $f(A) = A \backslash \{x\}$. Clearly, $f$ is
one-to-one. Suppose $f$ is onto. Then, since $Y \in \mathcal{I}$,
there exists $Z \in \mathcal{H}(x)$ such that $f(Z) = Y$. So $Z =
Y \cup \{x\} \in \mathcal{H}$, a contradiction since $Y$ is a
base. So $f$ is not onto. So the domain $\mathcal{H}(x)$ of $f$ is
smaller than the co-domain $\mathcal{I}$ of $f$. Since $|\mathcal{H}| = |\mathcal{H}(x)| + |\mathcal{I}|$, it follows that $|\mathcal{H}| > 2|\mathcal{H}(x)|$.~\hfill{$\Box$}\\
\\
%
\textbf{Proof of Theorem~\ref{main}.} We prove the result by induction on $m + n$. The basis is $m + n = 2$ with $m = n = 1$, in which case the result is trivial. Now consider $m + n > 2$. We may assume that $m \leq n$. If $m = 1$, then the result is trivial too, so we consider $m \geq 2$. If at least one of $\mathcal{G}$ and $\mathcal{H}$ is $\emptyset$ or $\{\emptyset\}$, then we trivially
have $|\mathcal{A}||\mathcal{B}| = 0 = |\mathcal{G}(1)||\mathcal{H}(1)|$. Thus we will assume that $\mathcal{G} \neq \emptyset$, $\mathcal{G} \neq \{\emptyset\}$, $\mathcal{H} \neq \emptyset$ and $\mathcal{H} \neq
\{\emptyset\}$. So each of $\mathcal{G}$ and $\mathcal{H}$
has at least one non-empty set.

As explained in Section~\ref{Compsection}, we apply left-compressions to $\mathcal{A}$ and $\mathcal{B}$ simultaneously until we obtain two compressed cross-intersecting families $\mathcal{A}^*$ and $\mathcal{B}^*$ such that $|\mathcal{A}^*| = |\mathcal{A}|$ and $|\mathcal{B}^*| = |\mathcal{B}|$. 
Since $\mathcal{G}$ and $\mathcal{H}$ are compressed, we have $\mathcal{A}^* \subseteq \mathcal{G}$ and $\mathcal{B}^* \subseteq \mathcal{H}$. We may therefore assume that
$\mathcal{A}$ and $\mathcal{B}$ are compressed.

Define $\mathcal{H}_0 = \{H \in \mathcal{H} \colon n \notin H\}$
and $\mathcal{H}_1 = \{H \backslash \{n\} \colon n \in H \in
\mathcal{H}\}$. Define $\mathcal{G}_0$, $\mathcal{G}_1$,
$\mathcal{A}_0$, $\mathcal{A}_1$, $\mathcal{B}_0$ and
$\mathcal{B}_1$ similarly. Since $\mathcal{A}$, $\mathcal{B}$,
$\mathcal{G}$ and $\mathcal{H}$ are compressed, we clearly have
that $\mathcal{A}_0$, $\mathcal{A}_1$, $\mathcal{B}_0$,
$\mathcal{B}_1$, $\mathcal{G}_0$, $\mathcal{G}_1$,
$\mathcal{H}_0$ and $\mathcal{H}_1$ are compressed. Since
$\mathcal{G}$ and $\mathcal{H}$ are hereditary, we clearly have
that $\mathcal{G}_0$, $\mathcal{G}_1$, $\mathcal{H}_0$ and
$\mathcal{H}_1$ are hereditary. Obviously, we also have
$\mathcal{A}_0 \subseteq \mathcal{G}_0 \subseteq 2^{[m-1]}$,
$\mathcal{A}_1 \subseteq \mathcal{G}_1 \subseteq 2^{[m-1]}$,
$\mathcal{B}_0 \subseteq \mathcal{H}_0 \subseteq 2^{[n-1]}$ and
$\mathcal{B}_1 \subseteq \mathcal{H}_1 \subseteq 2^{[n-1]}$. 

We have
\begin{equation} |\mathcal{B}| = |\mathcal{B}_0| +
|\mathcal{B}(n)| = |\mathcal{B}_0| + |\mathcal{B}_1|. \label{0.1}
\end{equation}
Along the same lines, 
\begin{align} |\mathcal{H}(1)| &= |\mathcal{H}_0(1)| + |\{H \in \mathcal{H} \colon 1, n \in H\}| = |\mathcal{H}_0(1)| + |\{H \backslash \{n\} \colon 1, n \in H \in \mathcal{H}\}| \nonumber \\
&= |\mathcal{H}_0(1)| + |\mathcal{H}_1(1)|. \label{0.2}
\end{align}

Suppose $m < n$. Clearly, $\mathcal{A}$ and $\mathcal{B}_0$ are
cross-intersecting. Since $m < n$, no set in $\mathcal{A}$ contains $n$, and hence $\mathcal{A}$ and $\mathcal{B}_1$ are cross-intersecting too. Thus, by the induction hypothesis,
$|\mathcal{A}||\mathcal{B}_q| \leq
|\mathcal{G}(1)||\mathcal{H}_q(1)|$ for each $q \in \{0, 1\}$. Together with (\ref{0.1}) and (\ref{0.2}), this gives us
\begin{align} |\mathcal{A}||\mathcal{B}| &= |\mathcal{A}|(|\mathcal{B}_0| + |\mathcal{B}_1|) = |\mathcal{A}||\mathcal{B}_0| + |\mathcal{A}||\mathcal{B}_1| \nonumber \\
&\leq |\mathcal{G}(1)||\mathcal{H}_0(1)| +
|\mathcal{G}(1)||\mathcal{H}_1(1)| = |\mathcal{G}(1)||\mathcal{H}(1)|, \nonumber
\end{align}
as required.

Now suppose $m=n$. Similarly to (\ref{0.1}) and (\ref{0.2}), we have

\begin{align} |\mathcal{A}| &= |\mathcal{A}_0| + |\mathcal{A}_1|, \label{0.3} \\
|\mathcal{G}(1)| &= |\mathcal{G}_0(1)| + |\mathcal{G}_1(1)|. \label{0.4}
\end{align}
Clearly, $\mathcal{A}_0$ and $\mathcal{B}_0$ are cross-intersecting, and since $n = m$, so are $\mathcal{A}_0$ and $\mathcal{B}_1$, and also $\mathcal{A}_1$ and $\mathcal{B}_0$.

Suppose $\mathcal{A}_1$ and $\mathcal{B}_1$ are
cross-intersecting too. Then, by the induction hypothesis,
$|\mathcal{A}_p||\mathcal{B}_q| \leq
|\mathcal{G}_p(1)||\mathcal{H}_q(1)|$ for any $p, q \in \{0, 1\}$. Together with (\ref{0.1})--(\ref{0.4}), this gives us 
\begin{align} |\mathcal{A}||\mathcal{B}| &=
(|\mathcal{A}_0| + |\mathcal{A}_1|)(|\mathcal{B}_0| +
|\mathcal{B}_1|) = |\mathcal{A}_0||\mathcal{B}_0| +
|\mathcal{A}_0||\mathcal{B}_1| + |\mathcal{A}_1||\mathcal{B}_0| +
|\mathcal{A}_1||\mathcal{B}_1| \nonumber \\
&\leq |\mathcal{G}_0(1)||\mathcal{H}_0(1)| +
|\mathcal{G}_0(1)||\mathcal{H}_1(1)| +
|\mathcal{G}_1(1)||\mathcal{H}_0(1)| +
|\mathcal{G}_1(1)||\mathcal{H}_1(1)|\nonumber \\
&= (|\mathcal{G}_0(1)| + |\mathcal{G}_1(1)|)(|\mathcal{H}_0(1)| +
|\mathcal{H}_1(1)|) = |\mathcal{G}(1)||\mathcal{H}(1)|, \label{1}
\end{align}
as required.

Now suppose $\mathcal{A}_1$ and $\mathcal{B}_1$ are not cross-intersecting. Then there exists $A \in \mathcal{A}_1$
such that $A \cap B = \emptyset$ for some $B \in \mathcal{B}_1$.
Let $\mathcal{C} = \{A \in \mathcal{A}_1 \colon A \cap B =
\emptyset \mbox{ for some } B \in \mathcal{B}_1\}$, and let $A_1,
\dots, A_k$ be the distinct sets in $\mathcal{C}$. For each $i
\in [k]$, let $B_i = [n-1] \backslash A_i$, $A_i' = A_i \cup \{n\}$, $B_i' = B_i \cup \{n\}$. For each $i \in [k]$, $A_i' \in \mathcal{A}$ since $A_i \in \mathcal{A}_1$.

Let $j \in [k]$. Let $\mathcal{D}_j = \{B \in \mathcal{B}_1 \colon A_j \cap B = \emptyset\}$. Suppose there exists $B \in \mathcal{D}_j$ such that $B \neq B_j$. Then $B \subsetneq [n-1] \backslash A_j$ and hence $[n-1] \backslash (A_j \cup B) \neq \emptyset$. Let $c \in [n-1] \backslash (A_j \cup B)$. Since $B \in \mathcal{B}_1$, $B \cup \{n\} \in \mathcal{B}$. Let $C = \delta_{c,n}(B \cup \{n\})$. Since $c \notin B \cup \{n\}$, $C = B \cup \{c\}$. Since $\mathcal{B}$ is compressed, $C \in \mathcal{B}$. However, since $c \notin A_j'$ and $A_j \cap B = \emptyset$, we have $A_j' \cap C = \emptyset$, which is a contradiction as $\mathcal{A}$ and $\mathcal{B}$ are cross-intersecting. So $\mathcal{D}_j \subseteq \{B_j\}$. By definition of $A_j$, $\mathcal{D}_j \neq \emptyset$ and hence $\mathcal{D}_j = \{B_j\}$.

We have therefore shown that
\begin{equation}\mbox{for each $i \in [k]$, $B_i$ is the unique set in $\mathcal{B}_1$ that does not intersect $A_i$,} \label{main8.1}
\end{equation}
and $B_i' \in \mathcal{B}$ (since $B_i \in \mathcal{B}_1$).
It follows from (\ref{main8.1}) that
\begin{equation}\mbox{for each $i \in [k]$, $A_i$ is the unique set in $\mathcal{A}_1$ that does not intersect $B_i$.} \label{main8.2}
\end{equation}
Indeed, suppose $A \in \mathcal{A}_1$ and $j \in [k]$ such that $A \cap B_j = \emptyset$. Then $A \in \mathcal{C}$. So $A = A_i$ for some $i \in [k]$. By (\ref{main8.1}), $B_j = B_i$. So $i = j$ as $B_1, \dots, B_k$ are distinct. So $A = A_j$.

Since $\mathcal{A}$ and $\mathcal{B}$ are compressed,
\begin{equation} \mbox{ for any $h \in [n-1]$ and any $i \in [k]$,  $\delta_{h,n}(A_i') \in \mathcal{A}$ and $\delta_{h,n}(B_i') \in \mathcal{B}$.} \label{main9}
\end{equation}

Let $r = \lfloor k/2 \rfloor$. Let
\[\begin{array}{ll}
\mathcal{A}_0' = \mathcal{A}_0 \cup \{A_i \colon i \in [r]\},
\quad \quad \quad & \mathcal{A}_1' = \mathcal{A}_1 \backslash
\{A_j \colon j \in [k] \backslash [r]\},\\
\mathcal{B}_0' = \mathcal{B}_0 \cup \{B_j \colon j \in [k]
\backslash [r]\}, & \mathcal{B}_1' = \mathcal{B}_1 \backslash
\{B_i \colon i \in [r]\}.
\end{array} \]
By (\ref{main8.1}) and (\ref{main8.2}), for any $p, q \in \{0, 1\}$, $\mathcal{A}_p'$ and $\mathcal{B}_q'$ are cross-intersecting. Since $\mathcal{G}$ and $\mathcal{H}$ are hereditary, we have $A_1, \dots, A_k \in \mathcal{G}_0$ and $B_1, \dots, B_k \in \mathcal{H}_0$, and hence $\mathcal{A}_0' \subseteq \mathcal{G}_0$ and $\mathcal{B}_0'
\subseteq \mathcal{H}_0$. Obviously, $\mathcal{A}_1' \subseteq
\mathcal{G}_1$ and $\mathcal{B}_1' \subseteq \mathcal{H}_1$.
Thus, by the induction hypothesis,
\begin{equation} |\mathcal{A}_p'||\mathcal{B}_q'| \leq
|\mathcal{G}_p(1)||\mathcal{H}_q(1)| \quad \mbox{for any } p, q
\in \{0,1\}.\label{1.5}
\end{equation}
As in the calculation in (\ref{1}), this gives us
\begin{equation}
(|\mathcal{A}_0'| + |\mathcal{A}_1'|)(|\mathcal{B}_0'| +
|\mathcal{B}_1'|) \leq |\mathcal{G}(1)||\mathcal{H}(1)| \label{2}
\end{equation}

Let $j \in [k]$. Since $B_j' \in \mathcal{B}(n)$ and
$A_j \cap B_j' = \emptyset$, the cross-intersection
condition gives us that $A_j \notin \mathcal{A}_0$. Similarly, $B_j \notin \mathcal{B}_0$. So
\begin{equation}A_1, \dots, A_k \notin \mathcal{A}_0 \quad \mbox{ and } \quad B_1, \dots, B_k \notin \mathcal{B}_0. \label{2.5}
\end{equation}
Therefore,
$|\mathcal{A}_0'| = |\mathcal{A}_0| + r$, $|\mathcal{A}_1'| =
|\mathcal{A}_1| + r - k$, $|\mathcal{B}_0'| = |\mathcal{B}_0| + k
- r$ and $|\mathcal{B}_1'| = |\mathcal{B}_1| - r$. So we have
\begin{align} |\mathcal{A}| &= |\mathcal{A}_0| + |\mathcal{A}_1| =
|\mathcal{A}_0'|+ |\mathcal{A}_1'| + k - 2r, \label{3}\\
|\mathcal{B}| &= |\mathcal{B}_0| + |\mathcal{B}_1| =
|\mathcal{B}_0'| + |\mathcal{B}_1'| + 2r - k. \label{4}
\end{align}

Suppose $k$ is even. Then $k = 2r$. By (\ref{3}) and (\ref{4}),
$|\mathcal{A}| = |\mathcal{A}_0'| + |\mathcal{A}_1'|$ and
$|\mathcal{B}| = |\mathcal{B}_0'| + |\mathcal{B}_1'|$. Thus, by (\ref{2}), $|\mathcal{A}||\mathcal{B}| \leq
|\mathcal{G}(1)||\mathcal{H}(1)|$, as required.

Suppose $k$ is odd. Then $k = 2r + 1$. By (\ref{3}) and
(\ref{4}), $|\mathcal{A}| = |\mathcal{A}_0'| + |\mathcal{A}_1'| +
1$ and $|\mathcal{B}| = |\mathcal{B}_0'| + |\mathcal{B}_1'| - 1$.
Let 
%
%
\begin{align} \mathcal{A}_0'' &= \mathcal{A}_0' \cup \{A_{r+1}\} = \mathcal{A}_0 \cup \{A_i \colon i \in [r+1]\}, \nonumber \\
\mathcal{A}_1'' &= \mathcal{A}_1' \cup \{A_{r+1}\} = \mathcal{A}_1 \backslash \{A_j \colon j \in [k] \backslash [r+1]\}, \nonumber \\
\mathcal{B}_0'' &= \mathcal{B}_0' \backslash \{B_{r+1}\} = \mathcal{B}_0 \cup \{B_j \colon j \in [k]
\backslash [r+1]\}, \nonumber \\
\mathcal{B}_1'' &= \mathcal{B}_1' \backslash \{B_{r+1}\} = \mathcal{B}_1 \backslash \{B_i \colon i \in [r+1]\}. \nonumber
\end{align}
Similarly to $\mathcal{A}_0'$, $\mathcal{A}_1'$, $\mathcal{B}_0'$ and $\mathcal{B}_1'$, we have that for any $p, q \in \{0,1\}$, $\mathcal{A}_p''$ and $\mathcal{B}_q''$ are cross-intersecting, $\mathcal{A}_p'' \subseteq \mathcal{G}_p$ and $\mathcal{B}_q'' \subseteq \mathcal{H}_q$. Thus, by the induction hypothesis,
\begin{equation} |\mathcal{A}_p''||\mathcal{B}_q''| \leq
|\mathcal{G}_p(1)||\mathcal{H}_q(1)| \quad \mbox{for any } p, q
\in \{0,1\}. \label{5}
\end{equation}
For each $p \in \{0,1\}$, let $a_p = |\mathcal{A}_p|$, $a_p' =
|\mathcal{A}_p'|$, $a_p'' = |\mathcal{A}_p''|$, $b_p =
|\mathcal{B}_p|$, $b_p' = |\mathcal{B}_p'|$, $b_p'' =
|\mathcal{B}_p''|$, $g_p = |\mathcal{G}_p(1)|$, $h_p =
|\mathcal{H}_p(1)|$. Recall that $A_1, \dots, A_k \in \mathcal{A}_1 \backslash \mathcal{A}_0$ (by (\ref{2.5})) and $B_1, \dots, B_k \in \mathcal{B}_1 \backslash \mathcal{B}_0$ (by (\ref{main8.1}) and (\ref{2.5})), and so we have $a_0'' = a_0' + 1$, $a_1'' = a_1' + 1$, $b_0'' = b_0' - 1$ and $b_1'' = b_1' - 1$. Together with (\ref{1.5}) and (\ref{5}), this gives us
\begin{align} a_0'b_0' &\leq g_0h_0, \label{6.1} \\
(a_0'+1)(b_0'-1) = a_0''b_0'' &\leq g_0h_0, \label{6.2} \\[2mm]
a_0'b_1' &\leq g_0h_1, \label{6.3} \\
(a_0'+1)(b_1'-1) = a_0''b_1'' &\leq g_0h_1, \label{6.4} \\[2mm]
a_1'b_0' &\leq g_1h_0, \label{6.5} \\
(a_1'+1)(b_0'-1) = a_1''b_0'' &\leq g_1h_0, \label{6.6} \\[2mm]
a_1'b_1' &\leq g_1h_1, \label{6.7} \\
(a_1'+1)(b_1'-1) = a_1''b_1'' &\leq g_1h_1. \label{6.8}
\end{align}

Suppose $|\mathcal{A}| \geq |\mathcal{B}| + 1$. We have
\begin{align} |\mathcal{A}||\mathcal{B}| &= (a_0' + a_1' +
1)(b_0' + b_1' - 1) \quad \quad \mbox{(by (\ref{3}),
(\ref{4}))}\nonumber \\
&= (a_0' + a_1')(b_0' + b_1') - (a_0' + a_1') + (b_0' + b_1') - 1
\nonumber \\
&\leq |\mathcal{G}(1)||\mathcal{H}(1)| - (|\mathcal{A}| - 1) +
(|\mathcal{B}| + 1) - 1 \quad \quad \mbox{(by (\ref{2}),
(\ref{3}), (\ref{4}))} \nonumber \\
&\leq |\mathcal{G}(1)||\mathcal{H}(1)|. \nonumber
\end{align}

Next, suppose $|\mathcal{A}| \leq |\mathcal{B}| - 1$. We have
\begin{align} |\mathcal{A}||\mathcal{B}| &= (a_0'' + a_1'' -
1)(b_0'' + b_1'' + 1) \quad \quad \mbox{(by (\ref{3}),
(\ref{4}))}\nonumber \\
&= (a_0'' + a_1'')(b_0'' + b_1'') + (a_0'' + a_1'') - (b_0'' +
b_1'') - 1 \nonumber \\
&= (a_0''b_0'' + a_0''b_1'' + a_1''b_0'' + a_1''b_1'') +
(|\mathcal{A}| + 1) - (|\mathcal{B}| - 1) - 1 \quad \quad
\mbox{(by (\ref{3}), (\ref{4}))} \nonumber \\
&\leq (g_0h_0 + g_0h_1 + g_1h_0 + g_1h_1) + |\mathcal{A}| + 1 -
|\mathcal{B}|  \quad \quad  \mbox{(by (\ref{6.2}),
(\ref{6.4}), (\ref{6.6}), (\ref{6.8}))} \nonumber\\
&= |\mathcal{G}(1)||\mathcal{H}(1)| + |\mathcal{A}| + 1 -
|\mathcal{B}| \quad \quad \mbox{(as in (\ref{1}))} \nonumber \\
&\leq |\mathcal{G}(1)||\mathcal{H}(1)|. \nonumber
\end{align}

Finally, suppose $|\mathcal{A}| = |\mathcal{B}|$. Then, by
(\ref{3}) and (\ref{4}), $a_0' + a_1' + 1 = b_0' + b_1' -
1$ and hence
\begin{equation} a_0' + a_1' = b_0' + b_1' - 2. \label{7}
\end{equation}
Also by (\ref{3}) and (\ref{4}), we have $(a_0' + a_1')(b_0' +
b_1') = (|\mathcal{A}| - 1)(|\mathcal{B}| + 1) = (|\mathcal{A}| -
1)(|\mathcal{A}| + 1) = |\mathcal{A}|^2 - 1 =
|\mathcal{A}||\mathcal{B}| - 1$ and hence
\begin{equation} a_0'b_0' + a_0'b_1' + a_1'b_0' + a_1'b_1' =
|\mathcal{A}||\mathcal{B}| - 1. \label{8}
\end{equation}
As in (\ref{1}), we have
\begin{equation} g_0h_0 + g_0h_1 + g_1h_0 + g_1h_1 =
|\mathcal{G}(1)||\mathcal{H}(1)|. \label{9}
\end{equation}
We now split the problem into the following cases.\medskip
\\
\textit{Case 1: $a_0' \neq a_1'$}. So either $a_0' \geq a_1' + 1$ or
$a_0' \leq a_1' - 1$.

Suppose $a_0' \geq a_1' + 1$. Then (\ref{6.1}),
(\ref{6.3}), (\ref{6.6}) and (\ref{6.8}) give us
\begin{align} & a_0'b_0' + a_0'b_1' + (a_1'+1)(b_0'-1) +
(a_1'+1)(b_1'-1) \leq g_0h_0 + g_0h_1 + g_1h_0 + g_1h_1 \nonumber
\\
\Rightarrow \; & (a_0'b_0' + a_0'b_1' + a_1'b_0' + a_1'b_1') +
(b_0' + b_1' - 2) - 2a_1' \leq |\mathcal{G}(1)||\mathcal{H}(1)|
\quad \quad \mbox{(by (\ref{9}))} \nonumber \\
\Rightarrow \; & (|\mathcal{A}||\mathcal{B}| - 1) + (a_0' + a_1')
- 2a_1' \leq |\mathcal{G}(1)||\mathcal{H}(1)| \quad \quad
\mbox{(by (\ref{7}), (\ref{8}))} \nonumber \\
\Rightarrow \; & |\mathcal{A}||\mathcal{B}| \leq
|\mathcal{G}(1)||\mathcal{H}(1)| \quad \quad \mbox{(as $a_0' \geq
a_1' + 1$)}.\nonumber
\end{align}

Now suppose $a_0' \leq a_1' - 1$. Then (\ref{6.2}),
(\ref{6.4}), (\ref{6.5}) and (\ref{6.7}) give us
\begin{align} & (a_0'+1)(b_0'-1) + (a_0'+1)(b_1'-1) + a_1'b_0' +
a_1'b_1' \leq g_0h_0 + g_0h_1 + g_1h_0 + g_1h_1 \nonumber \\
\Rightarrow \; & (a_0'b_0' + a_0'b_1' + a_1'b_0' + a_1'b_1') +
(b_0' + b_1' - 2) - 2a_0' \leq |\mathcal{G}(1)||\mathcal{H}(1)|
\quad \quad \mbox{(by (\ref{9}))} \nonumber \\
\Rightarrow \; & (|\mathcal{A}||\mathcal{B}| - 1) + (a_0' + a_1')
- 2a_0' \leq |\mathcal{G}(1)||\mathcal{H}(1)| \quad \quad
\mbox{(by (\ref{7}), (\ref{8}))} \nonumber \\
\Rightarrow \; & |\mathcal{A}||\mathcal{B}| \leq
|\mathcal{G}(1)||\mathcal{H}(1)| \quad \quad \mbox{(as $a_0' \leq
a_1' - 1$).} \nonumber
\end{align}
\textit{Case 2: $b_0' \neq b_1'$}. So either $b_0' \geq b_1' + 1$ or
$b_0' \leq b_1' - 1$.

Suppose $b_0' \geq b_1' + 1$. Then (\ref{6.2}),
(\ref{6.3}), (\ref{6.6}) and (\ref{6.7}) give us
\begin{align} & (a_0'+1)(b_0'-1) + a_0'b_1' + (a_1'+1)(b_0'-1) +
a_1'b_1' \leq g_0h_0 + g_0h_1 + g_1h_0 + g_1h_1 \nonumber \\
\Rightarrow \; & (a_0'b_0' + a_0'b_1' + a_1'b_0' + a_1'b_1') +
(b_0' - a_0' - a_1' - 2) + b_0' \leq
|\mathcal{G}(1)||\mathcal{H}(1)| \quad \quad \mbox{(by (\ref{9}))}
\nonumber \\
\Rightarrow \; & (|\mathcal{A}||\mathcal{B}| - 1) - b_1' + b_0'
\leq |\mathcal{G}(1)||\mathcal{H}(1)| \quad \quad
\mbox{(by (\ref{7}), (\ref{8}))} \nonumber \\
\Rightarrow \; & |\mathcal{A}||\mathcal{B}| \leq
|\mathcal{G}(1)||\mathcal{H}(1)| \quad \quad \mbox{(as $b_0' \geq
b_1' + 1$)}.\nonumber
\end{align}

Now suppose $b_0' \leq b_1' - 1$. Then (\ref{6.1}),
(\ref{6.4}), (\ref{6.5}) and (\ref{6.8}) give us
\begin{align} & a_0'b_0' + (a_0'+1)(b_1'-1) + a_1'b_0' +
(a_1'+1)(b_1'-1) \leq g_0h_0 + g_0h_1 + g_1h_0 + g_1h_1 \nonumber
\\
\Rightarrow \; & (a_0'b_0' + a_0'b_1' + a_1'b_0' + a_1'b_1') +
(b_1' - a_0' - a_1' - 2) + b_1' \leq
|\mathcal{G}(1)||\mathcal{H}(1)| \quad \quad \mbox{(by (\ref{9}))}
\nonumber \\
\Rightarrow \; & (|\mathcal{A}||\mathcal{B}| - 1) - b_0'+ b_1'
\leq |\mathcal{G}(1)||\mathcal{H}(1)| \quad \quad
\mbox{(by (\ref{7}), (\ref{8}))} \nonumber \\
\Rightarrow \; & |\mathcal{A}||\mathcal{B}| \leq
|\mathcal{G}(1)||\mathcal{H}(1)| \quad \quad \mbox{(as $b_0' \leq
b_1' - 1$)}.\nonumber
\end{align}
\textit{Case 3: $a_0' = a_1'$ and $b_0' = b_1'$}. If at least one
of the inequalities (\ref{6.1}), (\ref{6.3}), (\ref{6.5}) and
(\ref{6.7}) is not an equality, then we have
\begin{align} & a_0'b_0' + a_0'b_1' + a_1'b_0' + a_1'b_1' <
g_0h_0 + g_0h_1 + g_1h_0 + g_1h_1 \nonumber \\
\Rightarrow \; & |\mathcal{A}||\mathcal{B}| - 1 <
|\mathcal{G}(1)||\mathcal{H}(1)| \quad \quad
\mbox{(by (\ref{8}), (\ref{9}))} \nonumber \\
\Rightarrow \; & |\mathcal{A}||\mathcal{B}| \leq
|\mathcal{G}(1)||\mathcal{H}(1)|.\nonumber
\end{align}
Now suppose equality holds in each of (\ref{6.1}), (\ref{6.3}),
(\ref{6.5}) and (\ref{6.7}); that is, $a_0'b_0' = g_0h_0$,
$a_0'b_1' = g_0h_1$, $a_1'b_0' = g_1h_0$ and $a_1'b_1' = g_1h_1$.
Then, since $a_0' = a_1'$ and $b_0' = b_1'$, we have $g_0h_0 =
g_0h_1 = g_1h_0 = g_1h_1$.

We will now show that $g_0 > 0$. Recall that $\mathcal{G} \neq
\emptyset$ and $\mathcal{G} \neq \{\emptyset\}$; that is, there
exists $G \in \mathcal{G}$ such that $G \neq \emptyset$. Suppose
$1 \in G$. Since $\mathcal{G}$ is hereditary, $G \backslash
\{n\} \in \mathcal{G}$. So we have $G \backslash \{n\} \in
\mathcal{G}_0(1)$ and hence $g_0 > 0$. Now suppose $1 \notin
G$. Let $x$ be the largest integer in $G$. Then $1 \in
\delta_{1,x}(G)$, $n \notin \delta_{1,x}(G)$ (by choice of $x$)
and $\delta_{1,x}(G) \in \mathcal{G}$ (since $\mathcal{G}$ is
compressed). So we have $\delta_{1,x}(G) \in \mathcal{G}_0(1)$
and hence again $g_0 > 0$, as required.

Similarly, $h_0 > 0$. So $g_0h_0 > 0$. Together with $g_0h_0 =
g_0h_1 = g_1h_0 = g_1h_1$, this gives us $g_0 = g_1 > 0$ and $h_0 = h_1
> 0$.

Suppose $[n] \notin \mathcal{G}$. Let $\mathcal{G}_2 = \{G \in
\mathcal{G} \colon 1, n \in G\}$. Clearly, $\mathcal{G}_1(1) = \{G
\backslash \{n\} \colon G \in \mathcal{G}_2\}$ and hence $|\mathcal{G}_2| = g_1$. Since $g_1 > 0$, we have $\mathcal{G}_1(1) \neq \emptyset$ and hence $\mathcal{G}_2 \neq \emptyset$. Let $M$ be a base of $\mathcal{G}_2$. Since $[n] \notin \mathcal{G}$, $[n] \backslash M \neq \emptyset$. Let $y \in [n] \backslash M$, and let $M' = \delta_{y,n}(M)$. Since $1, n \in M$, $y \notin M$, and $\mathcal{G}$ is compressed, we have $M' \in \mathcal{G}_0(1)$. Let $N$ be a base of $\mathcal{G}$ such that $M' \subseteq N$. Then $n \notin N$ because otherwise we obtain $M \cup \{y\} \subsetneq N \in \mathcal{G}_2$ (a contradiction since $M$ is a base of $\mathcal{G}_2$). So $N \in \mathcal{G}_0(1)$. Let $\mathcal{F} = \{G \backslash \{1\} \colon G \in \mathcal{G}(1)\}$. So $|\mathcal{F}| = |\mathcal{G}(1)| = |\mathcal{G}_0(1)| + |\mathcal{G}_2| = g_0 + g_1 = 2g_1$. Since $\mathcal{G}$ is hereditary, $\mathcal{F}$ is hereditary. Let $X = M \backslash \{1\}$ and $Y = N \backslash \{1\}$. Since $M$ and $N$ are bases of $\mathcal{G}$ containing $1$, $X$ and $Y$ are bases of $\mathcal{F}$. Since $n \in X$ and $n \notin Y$, Lemma~\ref{her-lemma} gives us $|\mathcal{F}(n)| < |\mathcal{F}|/2 = g_1$, which is a contradiction because $|\mathcal{F}(n)| = |\{G \backslash \{1\} \colon 1, n \in G \in \mathcal{G}\}| = |\mathcal{G}_2| = g_1$.

Therefore, $[n] \in \mathcal{G}$. Similarly, $[n] \in \mathcal{H}$. Since $\mathcal{G}$ and $\mathcal{H}$ are hereditary, $2^{[n]}
\subseteq \mathcal{G}$ and $2^{[n]} \subseteq \mathcal{H}$. But
$\mathcal{G}, \mathcal{H} \subseteq 2^{[n]}$. So $\mathcal{G} =
\mathcal{H} = 2^{[n]}$. Since $\mathcal{A}$ and $\mathcal{B}$ are cross-intersecting, $[n] \backslash A \notin \mathcal{B}$ for any $A \in \mathcal{A}$. Thus, $\mathcal{B} \subseteq 2^{[n]} \backslash \{[n] \backslash A \colon A \in \mathcal{A}\}$ and hence $|\mathcal{B}| \leq 2^n - |\mathcal{A}|$. Since $0 \leq \left(|\mathcal{A}| - \frac{1}{2}2^n \right)^2 = |\mathcal{A}|^2 - |\mathcal{A}|2^n + \frac{1}{4}\left( 2^n \right)^2$, we have $|\mathcal{A}|(2^n - |\mathcal{A}|) \leq \left( \frac{1}{2}2^n \right)^2$ and hence $|\mathcal{A}||\mathcal{B}| \leq \left( 2^{n-1} \right)^2 = \left|\left\{A \in 2^{[n]} \colon 1 \in A \right\}\right|^2 = |\mathcal{G}(1)| |\mathcal{H}(1)|$.~\hfill{$\Box$}
\\

Note that Theorem~\ref{main} generalises to the following result in the same way Theorem~\ref{mainthm} generalises to Theorem~\ref{maincor2}.

\begin{theorem} \label{maincor} If $n_1, \dots, n_k \in \mathbb{N}$, $k \geq 2$, $\mathcal{A}_i \subseteq \mathcal{H}_i \subseteq 2^{[n_i]}$ for each $i \in [k]$, $\mathcal{H}_1, \dots, \mathcal{H}_k$ are hereditary and compressed, and $\mathcal{A}_1, \dots, \mathcal{A}_k$ are cross-intersecting, then
\[\prod_{i=1}^k |\mathcal{A}_i| \leq \prod_{i=1}^k \left| \mathcal{H}_i(1) \right|,\]
and equality holds if $\mathcal{A}_i = \mathcal{H}_i(1)$ for each $i \in [k]$.

\end{theorem}


\end{document}